\documentclass{article}

\usepackage{delarray,verbatim,enumerate,a4wide}
\usepackage{amsmath,amsthm,amstext,amsbsy,amssymb,amsfonts,amscd}
\usepackage[all]{xy}
\usepackage[french]{babel}
\usepackage[T1]{fontenc}
\usepackage[utf8]{inputenc}

\usepackage{scalerel}

\usepackage{upgreek,bbm}

\usepackage[small]{titlesec}

\usepackage{color}
\definecolor{vert}{rgb}{0.1,0.4,0.2}
\usepackage[colorlinks=true,linkcolor=blue,citecolor=vert]{hyperref}

\usepackage{calligra}
\DeclareFontShape{T1}{calligra}{m}{n}{<->s*[0.95]callig15}{}
\DeclareMathAlphabet{\mathscr}{T1}{calligra}{m}{n}

\newtheorem{Th}{Théorème}[]
\newtheorem{Lem}[Th]{Lemme}
\newtheorem{Prop}[Th]{Proposition}
\newtheorem{Cor}[Th]{Corollaire}

\newtheorem{DTh}[Th]{Définition \& Théorème}

\newtheorem*{Th*}{Théorème}
\newtheorem*{Cor*}{Corollaire}
\newtheorem*{Def*}{Définition}
\newtheorem*{ThA}{Théorème A}
\newtheorem*{ThB}{Théorème B}
\newtheorem*{ThC}{Théorème C}

\def\Preuve{\noindent {\it Preuve.~}}

\def\Remarque{\smallskip\noindent {\it Remarque.~}}
\def\Remarques{\smallskip\noindent {\it Remarques.~}}

		\def\QQ{\mathbb Q}	
	\def\ZZ{\mathbb Z}		
\def\F2{\mathbb{F}_2}	\def\Z2{\mathbb{Z}_2}		
\def\Zl{{\mathbb{Z}_\ell}} 	\def\Ql{\mathbb{Q}_\ell}		

 		\def\P{\mathcal  P}		\def\U{\mathcal  U}	\def\F{\mathcal  F}	
\def\J{\mathcal  J}  		\def\C{\mathcal  C}		\def\R{\mathcal  R}		
\def\Dl{\mathcal  D\ssi{\!}\ell} 	\def\Pl{\mathcal  P\ssi{\!}\ell}  	\def\Cl{\mathcal  C\!\ell}
\def\E{\mathcal  E}					\def\D{\mathcal  D}

		\def\p{{\mathfrak p}}						
		\def\l{{\mathfrak l}}

\def\wi{\widetilde}				
	\def\div{\operatorname{div}}
	\def\deg{\operatorname{deg}}		\def\Ind{\operatorname{Ind}}
\def\Gal{\operatorname{Gal}}			
\def\Ker{\operatorname{Ker}}				

\newcommand{\order}{\raise0.8pt \hbox{${\scriptstyle \#}$}}
\newcommand\scale[2]{\vstretch{#1}{\hstretch{#1}{#2}}}

\newcommand\si[1]{\scale{.7}{#1}}	\newcommand\ssi[1]{\scale{.5}{#1}}
\newcommand\ph{{\phantom{*}}}		\newcommand\cc{{\scale{.8}{\rm c}}}
\newcommand\ab{{\scale{.8}{\rm ab}}}	\newcommand\lc{{\scale{.8}{\rm lc}}}	\newcommand\nr{{\scale{.8}{\rm nr}}}

\def\%{{\scale{.8}{\infty}}}	\newcommand\res{{\scale{.8}{\rm res}}}			\newcommand\bp{{\scale{.8}{\rm bp}}}

\makeatletter
\newcommand*\wt[2][0.2ex]{%
        \begingroup
        \mathchoice{\wt@helper{#1}{#2}{\displaystyle}{\textfont}}
                   {\wt@helper{#1}{#2}{\textstyle}{\textfont}}
                   {\wt@helper{#1}{#2}{\scriptstyle}{\scriptfont}}
                   {\wt@helper{#1}{#2}{\scriptscriptstyle}{\scriptscriptfont}}%
        \endgroup
        #2%
}
\newcommand*\wt@helper[4]{%
        \def\currentfont{\the#41}%
        \def\currentskewchar{\char\the\skewchar\currentfont}%
        \setbox\tw@\hbox{\currentfont$#2$\currentskewchar}%
        \dimen@ii\wd\tw@
        \setbox\tw@\hbox{\currentfont$#2${}\currentskewchar}%
        \advance\dimen@ii-\wd\tw@
        \rlap{\raisebox{-#1}{$\m@th#3\kern\dimen@ii\widetilde{\phantom{#2}}$}}%
}
\makeatother

\def\wE{\,\wt[0.1ex]{\!\mathcal E}}		\def\wU{\wt[0.2ex]{\mathcal U}}	\def\wC{\wt[0.2ex]{\mathcal C}}
\def\wJ{\,\wt[0.2ex]{\!\mathcal J}}	\def\wCl{\wt[0.1ex]{\mathcal C\!\ell}} \def\wDl{\wt[0.2ex]{\mathcal D\ssi{\!}\ell}}
				\def\wnu{\wt[0.1ex]{\nu}}	\def\wPl{\wt[0.2ex]{\mathcal P\ssi{\!}\ell}}
\def\wCap{\wt[0.3ex]{C\!ap}}			
							
\def\wH{\,\wt[0.1ex]{\!H}}

\begin{document}

\title{\Large\bf Classes logarithmiques imaginaires des corps abéliens}
\author{ Jean-François {\sc Jaulent} }
\date{}

\maketitle

{\footnotesize
\noindent{\bf Résumé.}
Nous nous intéressons aux $\ell$-groupes de classes imaginaires, classiques ou logarithmiques, d'un corps abélien $K$, en relation avec la théorie d'Iwasawa. Nous en donnons d'abord une description des $\varphi$-composantes pour chaque caractère $\ell$-adique irréductible imaginaire du groupe $\Gal(K/\QQ)$ et nous généralisons dans ce cadre le critère de minimalité de Gold sur les invariants lambda des corps quadratiques. Nous étudions d'abord le cas semi-simple $\ell\nmid[K:\QQ]$, puis le cas général en connection avec les formules de transition de genre pour les corps surcirculaires.
}\smallskip

{\footnotesize
\noindent{\bf Abstract.}
We are interested in classical and logarithmic imaginary classes of abelian number fields in connection with Iwasawa theory. For any given odd prime $\ell$ and any imaginary abelian number field $K$, we compute the $\varphi$-components of the logarithmic $\ell$-class group of $K$ associated to imaginary irreducible $\ell$-adic characters $\varphi$ of $\Gal(K/\QQ)$ and we extend to this situation the Gold criterium on lambda invariants. We first study the semi-simple case $\ell\nmid[K/\QQ]$, thus the non semi-simple case in connection with transition formulas of genus for surcircular fields.
}\smallskip


{\small
\tableofcontents
}

\section*{Introduction}
\addcontentsline{toc}{section}{Introduction}

Le pro-$\ell$-groupe $\,\Cl^{\,\si{log}}_K$ des classes logarithmiques d'un corps de nombres $K$ a été introduits dans \cite{J28} en liaison avec l'étude des $\ell$-symboles de Hilbert et des noyaux sauvages de la $K$-théorie pour un nombre premier arbitraire fixé $\ell$. Concrètement, il se présente, tout comme le $\ell$-groupes de classes au sens habituel, comme quotient du $\Zl$-module $\Dl_K$ construit sur les idéaux premiers de $K$ par l'image $\Pl_K$ du tensorisé $\ell$-adique $\R_K=\Zl\otimes_\ZZ K^\times$ du groupe multiplicatif de $K$, à ceci près qu'aux places $\l\mid\ell$ il convient de remplacer les valuations habituelles $\nu_\l(\cdot)$ par les valuations logarithmiques $\wt\nu_\l(\cdot)$ définies à partir des logarithmes $\ell$-adiques des valeurs absolues $\ell$-adiques $|\cdot|_\l$.\smallskip

Ces groupes ont une interprétation naturelle en terme de théorie d'Iwasawa, comme quotient des genres du module complètement décomposé $\,\wC_K$ sur la $\Zl$-extension cyclotomique $K_\%$ de $K$; de sorte, par exemple, que la conjecture de Gross-Kuz'min (pour $\ell$ et pour $K$) équivaut à postuler la finitude du sous-groupe $\,\wCl_K$ des classes logarithmiques de degré nul (cf. \cite{J28,Gra2005,J31}). Et, pour $K$ totalement réel, la conjecture de Greenberg \cite{Grb1973,Grb1976}  sur la trivialité de l'invariant {\em lambda} d'Iwasawa standard revient à affirmer que $\,\wCl_K$ capitule dans $K_\%$ \cite{J60}. Notons que la transcription aux classes logarithmiques du classique théorème d'Artin-Furtwängler \cite{J54,J70} montre que cette capitulation a bien lieu dans la pro-$\ell$-extension abélienne {\em localement} cyclotomique maximale $K_\%^{\si{\mathrm loc}}$ de $K$.\smallskip

Nous résumons ci-après les trois principaux résultats de ce travail:

\subsection*{Le cas semi-simple}

Lorsque $K$ est un corps abélien imaginaire de degré $d$ étranger à $\ell$, l'algèbre de Galois $\ell$-adique $\Zl[\Delta]$ du groupe $\Delta=\Gal(K/\QQ]$ est un anneau semi-local, produit direct d'extensions non-ramifiées $\ZZ_\varphi$ de $\Zl$, indexées par les caractères $\ell$-adiques irréductibles $\varphi$ de $\Delta$; et les projecteurs correspondant à cette décomposition sont donnés par les idempotents primitifs associés:

\centerline{$e_\varphi=\frac{1}{d}\sum_{\sigma\in\Delta}\varphi(\sigma)\sigma^{\si{-1}}$.}

Les $\ell$-groupes de classes, tant ordinaires que logarithmiques, s'écrivent ainsi comme produits de leurs $\varphi$-composantes isotypiques  respectives. En particulier, leurs composantes imaginaires s'obtiennent en sommant sur les caractères irréductibles imaginaires, i.e. agissant non trivialement sur la conjugaison complexe $\bar\tau$. Il vient ainsi (cf. Th. \ref{CLog},ii):

\begin{ThA}
Pour $\varphi$ imaginaire irréductible représenté dans l'induit $\chi^\ph_\ell=\Ind_{\Delta_\ell}^\Delta 1_{\Delta_\ell}$ du caractère unité du sous-groupe de décomposition de $\ell$ dans $K$, l'ordre $\wt h_\varphi^{\,\deg\varphi}$ de la $\varphi$-composante $\,\wCl_K^{\,e_\varphi}$ du $\ell$-groupe des classes logarithmiques est donné par la formule (où $x \sim y$ signifie $x/y\in\Zl^\times$):\smallskip

\centerline{$\wt h_\varphi \sim  \frac{1}{\deg\l} \log_\ell(N_{\si{K_\l/\Ql}}(\eta_\varphi))$,}\smallskip

\noindent où $\l$ est l'une quelconque des places de $K$ au-dessus de $\ell$; $h_\varphi$ l'ordre dans  le $\ell$-groupe des classes $\,\Cl_K$ du diviseur isotypique $\l_\varphi=\l^{e_\varphi}$ construit sur $\l$; et $\eta_\varphi$ un générateur $\varphi$-isotypique de $\l_\varphi^{h_\varphi}$.
\end{ThA}

Soient alors $\,\C_K^{\,e_\varphi} = \varprojlim \,\Cl_{K_n}^{\,e_\varphi}$ et $\,\wC_K^{\,e_\varphi} = \varprojlim \,\wCl_{K_n}^{\,e_\varphi}$ les limites projectives des $\varphi$-composantes des $\ell$-groupes de classes ordinaires et logarithmiques dans la  $\Zl$-extension cyclotomique $K_\%=\bigcup_nK_n$ de $K$. 
Notant $\lambda_K^\varphi=\dim_{\ZZ_\varphi}\C_K^{\,e_\varphi}$ et $\wi\lambda_K^\varphi=\dim_{\ZZ_\varphi}\wC_K^{\,e_\varphi}$, nous obtenons le critère suivant (Th. \ref{CGold}):

\begin{ThB}
Pour $\varphi$ imaginaire irréductible représenté dans $\chi^\ph_\ell$ on a $\,\lambda_K^\varphi\geqslant 1$ et les équivalences:\smallskip

\centerline{$\lambda_K^\varphi =1 \quad\Leftrightarrow\quad  \wt\lambda_K^\varphi =0 \quad\Leftrightarrow\quad \,\wCl_K^{\,e_\varphi} = 1$.}
\end{ThB}

Et la $\varphi$-composante de la pro-$\ell$-extension de Bertrandias-Payan est alors très simple (cf. \S \ref{Schema}).

\subsection*{Le cas non semi-simple}

Si maintenant $L$ est un corps abélien imaginaire de degré $d\ell^m$ avec $m>0$ et $\ell\nmid 2d$, le groupe de Galois $\Gal(L/\QQ)$ s'écrit canoniquement comme produit direct de son $\ell$-sous-groupe de Sylow $G$ d'ordre $\ell^m$ et de son unique sous-groupe $\Delta$ d'ordre $d$. On peut donc définir les $\varphi$-composantes comme précédemment à l'aide des idempotents $e_\varphi$ et étudier la propagation de $K=L^G$ à $L$ du critère de trivialité précédent. Notant $L_\%=\bigcup_n L_n$ la $\Zl$-extension cyclotomique de $L$, puis définissant  $\lambda_L^\varphi=\dim_{\ZZ_\varphi}\C_L^{\,e_\varphi}$ et $\wi\lambda_L^\varphi=\dim_{\ZZ_\varphi}\wC_L^{\,e_\varphi}$ comme plus haut, nous obtenons (cf. Th. \ref{NSS}):

\begin{ThC}
Pour chaque caractère $\ell$-adique irréductible imaginaire $\varphi$, contenu dans l'induit $\chi^\ph_\ell=\Ind_{\Delta_\ell}^\Delta 1_{\Delta_\ell}$ du caractère unité du sous-groupe de décomposition de $\ell$ dans l'extension $K/\QQ$, la $\varphi$-composante du module d'Iwasawa $\,\C_L$ est un $\ZZ_\varphi$-module libre de dimension  $\lambda_L^\varphi \geqslant d_\ell$, où $d_\ell$ désigne l'indice de décomposition de $\ell$ dans $L_\%/L$. Et on a les équivalences:\smallskip
$$
\lambda_L^\varphi =d_\ell \quad\Leftrightarrow\quad  \wt\lambda_L^\varphi =0 \quad\Leftrightarrow\quad \left\{
\aligned
&\,\wCl_K^{\,e_\varphi} = 1 \quad \textrm{et, pour tout } p \textrm { vérifiant } \varphi \mid \chi_p^\ph,\\
&L/K \text{ logarithmiquement non-ramifiée en }  p\,;
\endaligned
\right.
$$
où $\chi^\ph_p=\Ind_{\Delta_p}^\Delta 1_{\Delta_p}$ est l'induit du caractère unité du sous-groupe de décomposition de $p$ dans $K/\QQ$.
\end{ThC}

\newpage
\section{Classes logarithmiques et classes ordinaires}

Le pro-$\ell$-groupe des classes logarithmiques $\,\Cl^{\si{\,log}}_K$ a été introduit dans \cite{J28} par analogie avec le groupe des classes au sens habituel en  envoyant le tensorisé $\R_K=\Zl\otimes_\ZZ K^\times$ dans le $\Zl$-module des diviseurs $\Dl_K=\oplus_{\p}\Zl\,\p$ construit sur les premiers de $K$ par la famille $(\wnu_\p)_\p$ obtenue en remplaçant les valuations habituelles $\nu_\l$ aux premiers $\l$ qui divisent $\ell$ par les valuations $\ell$-adiques $\wnu_\l$ données (à normalisation près) par les logarithmes des valeurs absolues $\ell$-adiques $\wnu_\l(\cdot)= \log_\ell(|\cdot|)/\log_\ell(1+\ell)$:\smallskip

\centerline{$1 \rightarrow \wE_K \rightarrow  \R_K \overset{\wi\nu}{\longrightarrow}  \Dl_K \rightarrow \Cl_K^{\si{\,log}} \rightarrow 1$.}\smallskip

\noindent Dans la suite exacte obtenue le noyau $\,\wE_K$ à gauche est ainsi le groupe des {\em unités logarithmiques}; le conoyau $\,\Cl_K^{\si{\,log}}$ à droite, celui des {\em classes logarithmiques}. Contrairement aux groupes de classes et d'unités au sens habituel, ce sont donc par construction des objets $\ell$-adiques. Leur importance provient de leur interprétation par la Théorie $\ell$-adique du corps de classes (cf. \cite{Gra2005,J28,J31}).\smallskip

Pour voir cela, introduisons pour chaque premier $\p$ de $K$  le compactifié $\ell$-adique du groupe $K_\p^\times$ défini par $\R_\p=\varprojlim K_\p^\times/K_\p^{\times\ell^n}$; et notons $\J_K= \prod_\p^\res \R_\p$ le $\ell$-adifié du groupe des idèles.\smallskip

 Du point de vue local, le noyau $\U_\p$ de $\nu_\p$ dans $\R_\p$ (autrement dit le sous-groupe des unités de $\R_\p$) est le groupe de normes associé à la $\Zl$-extension non ramifiée de $K_\p$; tandis que le noyau $\wU_\p$ de $\wnu_\p$ (i.e. le sous-groupe des unités logarithmiques locales) correspond, lui, à sa $\Zl$-extension cyclotomique. En d'autres termes, une $\ell$-extension abélienne de $K_\p$ est {\em logarithmiquement non-ramifiée} si et seulement si elle est contenue dans la $\Zl$-extension cyclotomique $K^c_\p$ de $K_\p$.\smallskip
 
Du point de vue global, le $\ell$-groupe des classes d'idéaux $\,\Cl_K$ s'interprète comme groupe de Galois de la $\ell$-extension abélienne non ramifiée maximale $K^{\nr}$ de $K$; et le $\ell$-groupe des classes logarithmiques $\,\Cl_K^{\si{\,log}}\simeq \J_K/\prod_\p\wU_\p\,\R_K$ comme groupe de Galois de sa pro-$\ell$-extension abélienne localement cyclotomique maximale $K^\lc$. Le corps $K^\lc$ est ainsi la plus grande pro-$\ell$-extension abélienne de $K$ qui est complètement décomposée (en toutes ses places)  au-dessus de la $\Zl$-extension cyclotomique $K^\cc$. Et comme $K^\lc$ contient $K^\cc$, le groupe  $\,\Cl_K^{\si{\,log}}$ n'est jamais fini.\smallskip

Il est donc naturel d'introduire le sous-groupe $\wJ_K$ de $\J_K$ associé à $K^\cc$, noyau du morphisme {\em degré} de $\J_K$ sur $\Gal(K^\cc/K)\simeq\Zl$ et de définir le {\em $\ell$-groupe des diviseurs logarithmiques de degré nul} par $\wDl_K=\wJ_K/\wU_K\R_K$; puis le {\em $\ell$-groupe des classes logarithmiques de degré nul} par:\smallskip

\centerline{$\,\wCl_K =  \wDl_K/\wPl_K \simeq \wJ_K/\prod_\p\wU_\p\,\R_K \simeq \Gal(K^\lc/K^\cc)$,}\smallskip

\noindent en notant $\wPl_K=\Pl_K\simeq\R_K/\wE_K$ le sous-groupe principal de $\Dl_K$, qui est contenu dans $\wDl_K$.\smallskip

Or, comme expliqué dans \cite{J28,J31}, la conjecture de Gross-Kuz'min postule précisément la finitude du groupe $\,\wCl_K$ et elle est en particulier satisfaite par tous les corps abéliens. Tout comme les $\ell$-groupes de classes au sens ordinaire $\Cl_K$, les $\ell$-groupes de classes logarithmiques (de degré nul) sont donc finis (au moins pour les corps abéliens), mais ils en diffèrent par la contribution des {\em classes sauvages}, i.e. construites sur les places au-dessus de $\ell$. Dans le contexte considéré ici, où le corps abélien $K$ est pris de degré étranger à $\ell$, ils ont, en revanche, même quotient modéré:

\begin{Lem}\label{Mod}
Soient $K$ un corps abélien de degré $d$ et $\ell$ un nombre premier ne divisant pas $d$. Écrivons $\,\Cl_K^{\,\si{[\ell]}}$ et $\,\wCl_K^{\,\si{[\ell]}}$ les sous-groupes respectifs des $\ell$-groupe des classes ordinaires $\,\Cl_K$ ou logarith\-miques $\,\wCl_K$ construits sur les classes des places au-dessus de $\ell$. Alors, pour chaque caractère $\ell$-adique irréductible $\varphi$ du groupe $\Delta=\Gal(K/\QQ)$, les $\varphi$-composantes des quotients modérés $\,\Cl_K/\Cl_K^{\,\si{[\ell]}}$ et $\,\wCl_K/\wCl_K^{\,\si{[\ell]}}$ sont $\ZZ_\varphi$-isomorphes, ce que nous écrivons:\smallskip

\centerline{$\big(\Cl_K^\ph/\Cl_K^{\,\si{[\ell]}}\big)^{e_\varphi} \simeq \big(\wCl_K^\ph/\wCl_K^{\,\si{[\ell]}}\big)^{e_\varphi}$.}
\end{Lem}

\Preuve On a directement $\big(\Cl_K^\ph/\Cl_K^{\,\si{[\ell]}}\big)^{e_\varphi} \simeq \big(\Cl_K^{\,\si{log}}/\Cl_K^{\,\si{log\,[\ell]}}\big)^{e_\varphi}$, puisque les valuations ordinaires et logarithmiques coïncident aux places modérées. D'où l'isomorphisme de $\ZZ_\varphi$-modules:\smallskip

\centerline{$\big(\Cl_K^\ph/\Cl_K^{\,\si{[\ell]}}\big)^{e_\varphi} \simeq \big(\Cl_K^{\,\si{log}}/\Cl_K^{\,\si{log\,[\ell]}}\big)^{e_\varphi} = \big(\wCl_K^\ph/\wCl_Ki^{\,\si{[\ell]}}\big)^{e_\varphi} $,}\smallskip

\noindent pour tout $\varphi \ne 1$, puisque la restriction de degré ne concerne que la 1-composante. Or, par ailleurs:\smallskip

\centerline{$\Cl_K^{\,e_{\si{1}}}\simeq\Cl_\QQ^\ph=1$; \quad et semblablement \quad $\wCl_K^{\,e_{\si{1}}}\simeq\wCl_\QQ^\ph=1$.}

\newpage
\section{Expression du nombre de classes logarithmiques}

Le calcul des groupes de classes logarithmiques est implanté dans le système {\sc pari/gp}, comme expliqué dans \cite{BJ2016}. Mais dans le cas particulier des corps abéliens imaginaires semi-simples considérés ici, il est possible d'accéder plus directement aux ordres des $\varphi$-composantes, comme découvert fortuitement par Gras en dressant des tables numériques pour les corps quadratiques.\smallskip

Soient, en effet, $K$ abélien imaginaire et $\ell$ premier ne divisant pas $d=[K/\QQ]$. Désignons par $\chi_\ell^\ph$ l'induit à $\Delta=\Gal(K/\QQ)$ du caractère unité du sous-groupe de décomposition $\Delta_\ell$ de $\ell$.

Le sous-groupe sauvage $\D_K^{[\ell]}$ du $\ell$-adifié $\D_K^\ph=\Zl\otimes_\ZZ D_K^\ph$ du groupe des diviseurs de $K$ est le $\Zl[\Delta]$-module, isomorphe à $\Zl[\Delta/\Delta_\ell]$, engendré par l'une quelconque des places $\l$ au-dessus de $\ell$. Et le même résultat vaut pour le sous-groupe sauvage $\Dl_K^{[\ell]}$ du $\ell$-groupe des diviseurs logarithmiques $\Dl_K$, à cette seule différence que ce dernier est noté additivement et non plus multiplicativement.\smallskip

Convenons alors d'écrire $x \sim y$ lorsque $x/y$ est une unité $\ell$-adique. Cela étant, nous avons:

\begin{Th}\label{CLog}
Soient $K$ un corps abélien imaginaire, $\ell$ un nombre premier ne divisant pas l'ordre du groupe $\Delta=\Gal(K/\QQ)$ et $\chi^\ph_\ell=\Ind_{\Delta_\ell}^\Delta 1_{\Delta_\ell}$ l'induit à $\Delta$ du caractère unité du sous-groupe de décomposition de $\ell$. Alors:\smallskip
\begin{itemize}
\item[(i)] Pour chaque caractère $\ell$-adique irréductible $\varphi$ de $\Delta$ qui n'est pas représenté dans $\chi_\ell^\ph$ les $\varphi$-composantes respectives $\,\wCl_K^{\,e_\varphi}$ du $\ell$-groupe des classes logarithmiques et $\,\Cl_K^{\,e_\varphi}$ du $\ell$-groupe des classes ordinaire sont $\ZZ_\varphi$-isomorphes et ont, en particulier, le même ordre:\smallskip

\centerline{$\wt h_\varphi^{\,\deg\varphi} = h_\varphi^{\,\deg\varphi}$.}\smallskip

\item[(ii)] Pour $\varphi$ irréductible imaginaire représenté dans $\chi_\ell^\ph$, l'ordre $\wt h_\varphi^{\,\deg\varphi}$  de $\,\wCl_K^{\,e_\varphi}$ est donné par:\smallskip

\centerline{$\wt h_\varphi \sim  \frac{1}{\deg\l} \log_\ell(N_{\si{K_\l/\Ql}}(\eta_\varphi))$,}\smallskip

\noindent où $\l$ est l'une quelconque des places de $K$ au-dessus de $\ell$; $h_\varphi$ l'ordre dans $\,\Cl_K$ du diviseur isotypique $\l_\varphi=\l^{e_\varphi}$ construit sur $\l$; et $\eta_\varphi$ un générateur $\varphi$-isotypique de $\l_\varphi^{h_\varphi}$.
\end{itemize}
\end{Th}

\Preuve Observons d'abord que, l'anneau $\ZZ_\varphi$ étant non-ramifié sur $\Zl$, tout $\ZZ_\varphi$-module cyclique est de la forme $\ZZ_\varphi/\ell^a\ZZ_\varphi$ donc isomorphe comme $\Zl$-module à $(\Zl/\ell^a\Zl)^{\deg\varphi}$. Ainsi, puisque tout $\ZZ_\varphi$-module fini est un produit direct de sous-modules cycliques, son ordre est la puissance $\deg\varphi$-ième d'une puissance de $\ell$. Écrivons donc $h_\varphi^{\deg\varphi}$ l'ordre de la $\varphi$-composante $\,\Cl_K^{\,e_\varphi}$ du $\ell$-groupe des classes et $\wt h_\varphi^{\deg\varphi}$
celui de $\,\wCl_K^{\,e_\varphi}$. Cela étant:

$(i)$ Si $\varphi$ n'est pas représenté dans $\chi_\ell^\ph$, les $\ZZ_\varphi$-modules sauvages  $\big(\Cl_K^{\si{[\ell]]}}\big)^{e_\varphi}$ \!et $\big( \wCl_\varphi^{\si{[\ell]]}}\big)^{e_\varphi}$ \!sont tous deux triviaux et $\,\Cl_K^{\,e_\varphi}$ comme $\,\wCl_K^{\,e_\varphi}$ se réduisent à leurs quotients modérés respectifs, lesquels sont isomorphes d'après le Lemme \ref{Mod}; d'où l'égalité des ordres: $\wt h_\varphi^{\,\deg\varphi} = h_\varphi^{\,\deg\varphi}$.\smallskip

$(ii)$ Si, en revanche, $\varphi$ est représenté dans $\chi_\ell^\ph$, faisons choix de l'une quelconque des places $\l$ au-dessus de $\ell$; écrivons $[\l_\varphi]$ la classe du diviseur $\l_\varphi=\l^{e_\varphi}$ dans $\,\Cl_K$ et $[\wt\l_\varphi]$ celle du diviseur logarithmique $\wt\l_\varphi=e_\varphi\l$ qui est bien dans $\,\wCl_K$, puisque $\varphi$, qui est imaginaire, est différent de 1; notons enfin $w_\varphi$ l'ordre de $[\l_\varphi]$ et $\wt w_\varphi$ celui de $[\wt\l_\varphi]$, de sorte que nous avons:\smallskip

\centerline{$\big(\Cl_K^{\si{[\ell]}}\big)^{e_\varphi} \simeq \ZZ_\varphi/w_\varphi\ZZ_\varphi$ \quad \& \quad $\big(\wCl_K^{\si{[\ell]}}\big)^{e_\varphi} \simeq \ZZ_\varphi/\wt w_\varphi\ZZ_\varphi$}\smallskip

Par construction, nous avons $\l_\varphi^{w_\varphi}=(\varepsilon_\varphi)$ pour un générateur $\varepsilon_\varphi$ de la $\varphi$-composante  du $\ell$-adifié $\,\E'_K=\Zl\otimes_\ZZ E'_K$ du groupe des $\ell$-unités de $K$ (aux racines $\ell$-ièmes de l'unité éventuelles près). Et $\wt w_\varphi$ est donc le plus petit entier tel que l'on ait: $\wt w_\varphi \wt\l_\varphi\in \Pl_P^{\,e_\varphi}=\ZZ_\varphi\wt\div(\varepsilon_\varphi)$.
Or, comparant alors les valuations logarithmiques des deux termes, nous avons:\smallskip

\centerline{$\wt\nu_\l(\wt w_\varphi \wt\l_\varphi)=\dfrac{1}{d}\,\wt w_\varphi \sim \wt w_\varphi$ \quad et\quad $\wt\nu_\l(\varepsilon_\varphi) =\dfrac{1}{\deg\l}\, \log_\ell(N_{\si{K_\l/\Ql}}(\varepsilon_\varphi))$.}

\noindent D'où la valeur de $\wt w_\varphi$:

\centerline{$\wt w_\varphi \sim \dfrac{\log_\ell(N_{\si{K_\l/\Ql}}(\varepsilon_\varphi))}{\log_\ell(1+\ell)} \sim \frac{1}{\ell}\,\log_\ell(N_{\si{K_\l/\Ql}}(\varepsilon_\varphi))$.}

\noindent Puis, par le Lemme \ref{Mod}:

\centerline{$\wt h_\varphi =h_\varphi \, \dfrac{\wt w_\varphi}{w_\varphi} 
\sim \dfrac{h_\varphi}{w_\varphi}\, \big( \frac{1}{\deg\l}\log_\ell(N_{\si{K_\l/\Ql}}(\varepsilon_\varphi))\big)
= \frac{1}{\deg\l} \log_\ell(N_{\si{K_\l/\Ql}}(\eta_\varphi))$.}

\newpage
\section{Suites exactes des classes ambiges}

Soit maintenant $L/K$ une $\ell$-extension cyclique de groupe $G$. La suite exacte des classes ambige de Chevalley \cite{Chv} pour les $\ell$-groupes de classes dans l'extension $L/K$ s'obtient habituellement en partant de la suite exacte courte $1 \to \P_L \to \D_L \to \Cl_L \to 1$ qui définit $\Cl_L$ comme quotient du tensorisé $\D_L=\Zl\otimes_\ZZ D_L$ du groupe des diviseurs de $L$ par son sous-groupe principal $\P_L=\Zl\otimes_\ZZ P_L$, puis en comparant par le lemme du serpent la suite de cohomologie obtenue avec la même suite écrite pour $\Cl_K$ (cf. e.g. \cite{Gra2017a} ou \cite{J18}, Ch. III.1.1), ce qui conduit à la suite exacte canonique:\smallskip

\centerline{$1 \to \C\! ap_{L/K} \to \P_L^{\,G}/\P_K^\ph \to \D_L^{\,G}/\D_K^\ph \to \Cl_L^{\,G}/j_{\si{L/K}}(\Cl_K^\ph) \to H^1(G,\P_L^\ph) \to 1$.}\smallskip

\noindent où $j_{\si{L/K}}$ désigne le morphisme d'extension attaché à $L/K$ et $\,\C\!ap_{L/K}=\Ker j_{\si{L/K}}$ la capitulation.

Le groupe $ \P_L^{\,G}/\P_K^\ph$ s’interprète en termes d'unités depuis Chevalley comme quotient $_{\si{N}}\E_L/\E_L^{\si{I_G}}=H^1(G,\E_N)$ du noyau de la norme $N=N_{\si{L/K}}$ par l'image de l'augmentation et le groupe $H^1(G,\P_L^\ph)$ comme quotient $(\E_K\cap\ N_{\si{L/K}}(\R_L)/N_{\si{L/K}}(\E_L)$ des unités qui sont normes modulo les normes d'unités.\smallskip

Dans le contexte de cet article et pour $L$ absolument abélien le groupe $\Delta=\Gal(K/\QQ)$ opère sur chacun des termes de la suite, ce qui permet  de la spécialiser en les composantes isotypiques:\smallskip

\centerline{$1 \to \C\! ap_{L/K}^{e_\varphi} \to H^1(G,\E_N^{e_\varphi}) \to (\D_L^{\,G}/\D_K^\ph)^{e_\varphi} \to (\Cl_L^{\,G}/j(\Cl_K^\ph))^{e_\varphi} \to ((\E_K\cap\ N(\R_L)/N(\E_L))^{e_\varphi} \to 1$.}\smallskip

\noindent Ainsi, prenant pour $L$ un étage $K_n$ de la $\Zl$-extension cyclotomique $K^\cc$ et $\varphi$ imaginaire, on obtient:\smallskip

\centerline{$ \C\! ap_{K_n/K}^{e_\varphi} = H^1(G,\E_{K_n}^{e_\varphi}) =1$ \quad \& \quad $(\D_{K_n}^{\,G}/\D_K^\ph)^{e_\varphi} = (\Cl_{K_n}^{\,G}/j(\Cl_K^\ph))^{e_\varphi}$,}\smallskip

\noindent les unités imaginaires se réduisant aux racines de l'unité,  qui sont cohomologiquement triviales.
En particulier les classes ambiges imaginaires sont alors les classes d'ambiges, c'est-à-dire les classes des diviseurs imaginaires construits sur les étendus des diviseurs de $K$ et ceux au-dessus de $\ell$.\medskip

Le même formalisme vaut {\em mutatis mutandis} pour les classes logarithmiques. D'où, dans ce cas:\smallskip

\centerline{$1 \to \C\! ap_{L/K}^{\,\si{log}} \to \Pl_L^{\,G}/\Pl_K^\ph \to \Dl_L^{\,G}/\Dl_K^\ph \to \Cl_L^{\,\si{log}\,G}/j^{\,\si{log}}_{\si{L/K}}(\Cl_K^{\,\si{log}}) \to H^1(G,\Pl_L^\ph) \to 1$.}\smallskip

\noindent où $j^{\,\si{log}}_{\si{L/K}}$ désigne le morphisme d'extension et $\,\C\!ap^{\,\si{log}}_{L/K}=\Ker j^{\,\si{log}}_{\si{L/K}}$ la capitulation logarithmique. 

Ici encore le  groupe $ \Pl_L^{\,G}/\Pl_K^\ph$ s'interprète en termes d'unités logarithmiques comme quotient $_N\wE_L/\wE_L^{I_G}=H^1(G,\wE_N)$ et le groupe $H^1(G,\Pl_L^\ph)$ comme quotient $(\wE_K\cap\ N_{\si{L/K}}(\R_L)/N_{\si{L/K}}(\wE_L)$ des unités logarithmiques qui sont normes modulo les normes d'unités logarithmiques (cf. \cite{J28}). Mais une difficulté apparaît: si la capitulation logarithmique concerne évidemment les seules classes d'ordre fini et donc de degré nul, il n'est pas possible de remplacer directement le groupe $\Cl^{\,\si{log}}_L$ par son sous-groupe $\wCl_L$ car on n'a pas nécessairement $H^1(G,\wDl_L)=1$ bien qu'on ait $H^1(G,\Dl_L)=1$.\smallskip

En revanche, dans le contexte cyclotomique considéré plus haut, i.e. pour $L=K_n \subset K^\cc$, la restriction aux $\varphi$-composantes imaginaires fait disparaître la complication liée au degré (car celle-ci est spécifique à la 1-composante) de sorte qu'on a, comme précédemment:\smallskip

\centerline{$\wt{\C\! ap}_{K_n/K}^{e_\varphi} = H^1(G,\wE_{K_n}^{e_\varphi}) =1$ \quad \& \quad $(\wDl_{K_n}^{\,G}/\wDl_K^\ph)^{e_\varphi} = (\wCl_{K_n}^{\,G}/j(\Cl_K^\ph))^{e_\varphi}$,}\smallskip

\noindent avec ici $\,\wDl_{K_n}^{\,G}\!=\wDl_K^\ph$, puisque l'extension cyclotomique $K_n/K$ est logarithmiquement non-ramifiée.\medskip

En résumé, il vient:

\begin{Prop}
Soient $K$ un corps abélien imaginaire, $\ell$ un nombre premier ne divisant pas $d=[K:\QQ]$ et $K_n$ un étage fini de la $\Zl$-extension cyclotomique de $K$. Alors, pour chaque caractère $\ell$-adique irréductible imaginaire $\varphi$ du groupe $\Delta=\Gal(K/\QQ)$:\smallskip
\begin{itemize}
\item[(i)] Le morphisme d'extension $j_{\si{K_n/K}}$ envoie injectivement la $\varphi$-composante $\,\Cl_K^{\,e_\varphi}$ du $\ell$-groupe des classes de $K$ dans la $\varphi$-composante du sous-groupe ambige du $\ell$-groupe des classes de $K_n$; laquelle est $\ZZ_\varphi$-engendrée modulo $j_{\si{K_n/K}}(\Cl_K^{\,e_\varphi})$ par la classe $[\l^{e_\varphi}_n]$ du diviseur isotypique $\l^{e_\varphi}_n$ construit sur l'une quelconque des places $\l_n$ de $K_n$ au-dessus de $\ell$.\smallskip
\item[(ii)] Le morphisme d'extension logarithmique $\wt{\text \j} _{\si{K_n/K}}$ identifie la $\varphi$-composante $\,\wCl_K^{\,e_\varphi}$ du $\ell$-groupe des classes logarithmiques de $K$ à la $\varphi$-composante du sous-groupe ambige du $\ell$-groupe des classes logarithmiques de $K_n$.
\end{itemize}
\end{Prop}

\newpage
\section{Généralisation logarithmique du critère de Gold}

Le résultat qui suit peut être regardé comme le déploiement semi-simple de la généralisation par Gras \cite{Gra2024} du critère de Gold \cite{Gold1974} pour les corps quadratiques imaginaires.

\begin{Th}\label {CGold}
Soient $K$ un corps abélien imaginaire, $\ell$ un nombre premier ne divisant pas l'ordre du groupe $\Delta=\Gal(K/\QQ)$ et $\varphi$ un caractère $\ell$-adique irréductible imaginaire contenu dans l'induit $\chi^\ph_\ell=\Ind_{\Delta_\ell}^\Delta 1_{\Delta_\ell}$ du caractère unité du sous-groupe de décomposition de $\ell$.\smallskip

Alors l'invariant lambda d'Iwasawa du $\Lambda_\varphi$-module $\,\C_K^{\,e_\varphi} = \varprojlim \,\Cl_{K_n}^{\,e_\varphi}$ attaché à la $\Zl$-extension cyclotomique $K_\%=\bigcup_nK_n$ vérifie l'inégalité $\lambda_\varphi \geqslant 1$; et les assertions suivantes sont équivalentes:
\begin{itemize}\smallskip

\item[(i)] L'invariant structurel lambda du $\Lambda_\varphi$-module $\,\C_K^{\,e_\varphi} = \varprojlim \,\Cl_{K_n}^{\,e_\varphi}$ vaut: $\lambda_\varphi =1$.\smallskip

\item[(ii)] L'invariant structurel lambda du $\Lambda_\varphi$-module $\,\wC_K^{\,e_\varphi} = \varprojlim \,\wCl_{K_n}^{\,e_\varphi}$ vaut: $\wt\lambda_\varphi =0$.\smallskip

\item[(iii)] Le  $\Lambda_\varphi$-module $\,\wC_K^{\,e_\varphi}$ est pseudo-nul: $\,\wC_K^{\,e_\varphi} \sim 1$.\smallskip

\item[(iv)] Le  $\Lambda_\varphi$-module $\,\wC_K^{\,e_\varphi}$ est trivial: $\,\wC_K^{\,e_\varphi} = 1$.\smallskip

\item[(v)] La $\varphi$-composante du $\ell$-groupe des classes logarithmiques de $K$ est triviale: $\,\wCl_K^{\,e_\varphi} = 1$.
\end{itemize}
\end{Th}

\Preuve Vérifions tout d'abord l'inégalité $\lambda_\varphi \geqslant 1$: regardons pour cela le sous-module sauvage $\,\C_K^{\,\si{[\ell]}\,e_\varphi} = \varprojlim \,\Cl_{K_n}^{\,\si{[\ell]}\,e_\varphi}$ de $\,\C_K^{\,e_\varphi}$ construit sur les classes des idéaux au-dessus de $\ell$. Puisque, d'un côté, ceux-ci sont totalement ramifiés dans $K_\%/K$ et que, d'un autre côté, les $\ell$-unités imaginaires (sauf, peut-être, les racines $\ell^\%$-ièmes de l'unité) contenues dans $K_\%$ sont déjà dans $K$, la composante imaginaire de $\,\C_K$ s'identifie à celle du $\Zl[\Delta]$-module $\D_K^{\,\si{[\ell]}} \simeq\Zl[\Delta/\Delta_\ell]$ construit sur les idéaux de $K$ au-dessus de $\ell$. Il s'ensuit que chaque caractère $\ell$-adique irréductible imaginaire $\varphi$ est représenté exactement une fois dans  $\,\C_K^{\,\si{[\ell]}\,e_\varphi}$, donc au moins une fois dans $\,\C_K^{\,e_\varphi}$, comme annoncé.

Cela étant, nous avons successivement:
\begin{itemize}
\item[$(i) \Leftrightarrow (ii)$] en vertu du Lemme \ref{Mod}, puisque, les places logarithmiques étant inertes dans la tour $K_\%/K$, le groupe $\,\wC_K^{\,e_\varphi} = \varprojlim \,\wCl_{K_n}^{\,e_\varphi}$ s'identifie à son quotient modéré $\,\C'{}_{\!K}^{\,e_\varphi} = \varprojlim \,\wCl'{}_{\!K_n}^{\,e_\varphi}$ donc a même paramètre lambda: $\wt\lambda_\varphi=\lambda'_\varphi=\lambda_\varphi-1$;

\item[$(ii) \Leftrightarrow (iii)$] puisque le paramètre $\wt\mu_\varphi=\mu_\varphi$ est nul par le Théorème de Ferrero et Washington;

\item[$(iii) \Leftrightarrow (iv)$] en l'absence de capitulation imaginaire pour les classes logarithmiques (comme pour les classes ordinaires, faute d'unités imaginaires), puisque le sous-module fini de $\,\wC_K^{\,e_\varphi}$ s'identifie à la limite projective des sous-groupes de capitulation $\wCap_K^{\,e_\varphi}$;

\item[$(iv) \Leftrightarrow (v)$] puisque  $\,\wCl_K^{\,e_\varphi}={}^\Gamma\wC_K^{\,e_\varphi}$ n'est autre que le quotient des co-points fixes de $\,\wC_K^{\,e_\varphi}$ relativement au groupe procyclique $\Gamma=\Gal(K_\%/K)$.
\end{itemize}

\begin{Cor}
Lorsque les conditions équivalentes du Théorème sont vérifiées, pour tout $n \geqslant 0$ la $\varphi$-composante $\,\Cl_{K_n}^{\,e_\varphi}$ du $\ell$-groupe des classes de $K_n$ est $Z_\varphi$-monogène, engendrée par la classe du diviseur $\l_n^{e_\varphi}$ et d'ordre $\ell^{(n+\nu_\varphi)\deg\,\varphi}$ où $\ell^{\nu_\varphi}$ est l'ordre dans $\,\Cl_K$ des premiers $\l$ au-dessus de $\ell$.
\end{Cor}

\Preuve Pour $\lambda_\varphi=1$, l'ordre de la $\varphi$-composante de $\,\Cl_{K_n}$ est donnée asymptotiquement par:\smallskip

\centerline{$|\,\Cl_{K_n}^{\,e_\varphi}|=\ell^{(n+\nu_\varphi)\deg\varphi}$, pour un $\nu_\varphi$ a priori dans $\ZZ$.}\smallskip

\noindent Prenant $m> n \gg 0$ et appliquant à $\,\Cl_{K_m}^{\,e_\varphi}$ la formule des classes ambiges dans l'extension $K_m/K_n$, disons de groupe $\Gamma_{\!m:n}=\Gal(K_m/K_n)$, nous obtenons:\smallskip

\centerline{$|\,\Cl_{K_m}^{\,e_\varphi\,\Gamma_{\!m:n}}|=|\,\Cl_{K_n}^{\,e_\varphi}|\,\ell^{(m-n)\deg\varphi}=\ell^{(m+\nu_\varphi)\deg\varphi} =|\,\Cl_{K_m}^{\,e_\varphi}|$.}\smallskip

\noindent Ainsi $\,\Cl_{K_m}^{\,e_\varphi}$ coïncide avec son sous-groupe ambige, i.e. (en l'absence d'unités imaginaires) avec le sous-groupe des classes d'ambiges, qui est engendré conjointement par le sous-groupe $j_{\si{K_m/K_n}}(\Cl_{K_n}^{\,e_\varphi})$ formé des classes étendues de $K_n$ et le sous-groupe sauvage $\,\Cl_{K_m}^{\si{\,[\ell]}\,e_\varphi}$. Or, le premier est inutile par:\smallskip

\centerline{$j_{\si{K_m/K_n}}(\Cl_{K_n}^{\,e_\varphi})=j_{\si{K_m/K_n}}(N_{\si{K_m/K_n}}(\Cl_{K_m}^{\,e_\varphi}))=j_{\si{K_m/K_n}}(N_{\si{K_m/K_n}}(\Cl_{K_m}^{\,e_\varphi\,\Gamma_{\!m:n}}))=\Cl_{K_m}^{\,e_\varphi\,\ell^{m-n}} \!\!\subset \Cl_{K_m}^{\,e_\varphi\,\ell}$.}\medskip

\noindent Il suit $\,\Cl_{K_m}^{\,e_\varphi}= \Cl_{K_m}^{\,\si{[\ell]}\,e_\varphi}$ pour $m \gg 0$; puis $\,\Cl_{K_n}^{\,e_\varphi}=N_{\si{K_m/K_n}}(\Cl_{K_m}^{\,e_\varphi})=N_{\si{K_m/K_n}}( \Cl_{K_m}^{\,\si{[\ell]}\,e_\varphi}) = \,\Cl_{K_n}^{\,\si{[\ell]}\,e_\varphi}$ quel que soit $n\geqslant 0$. 
Ainsi $\,\Cl_{K_n}^{\,e_\varphi}$ est $Z_\varphi$-monogène, engendré par la classe de $\l^{e_\varphi}_n$ d'ordre $\ell^{n\deg\varphi} |\,\Cl^{\si{\,[\ell]}\,e_\varphi}_K|$.

\newpage
\section{Construction des $\ZZ_\varphi$-extensions absolument galoisiennes}

Pour les caractères réels, il n'est d'autre $\Zl$-extension isotypique que la cyclotomique. Mais:

\begin{Prop}
Soit $K$ un corps abélien imaginaire et $\ell$ un nombre premier ne divisant pas son degré.
Pour chaque caractère $\ell$-adique irréductible imaginaire $\varphi$ du groupe $\Delta=\Gal(K/\QQ)$ le corps $K$ possède une unique $\Zl$-extension multiple $Z^\varphi$ dont le groupe $\Gamma_{\!\varphi}=\Gal(Z^\varphi/K)$ est $\varphi$-isotypique.\par

Le groupe $\Gamma_{\!\varphi}=\gamma_\varphi ^{\ZZ_\varphi}$ est ainsi un $\ZZ_\varphi$-module libre de rang 1 sur lequel $\Delta$ agit par conjugaison:\smallskip

\centerline{$\tau\gamma_\varphi^z\tau^{\si{-1}}=\gamma_\varphi^{\rho_\varphi^\tau(z)},$ \quad pour $\tau\in\Delta$ et $z \in \ZZ_\varphi$.}\smallskip

Nous disons que $Z^\varphi$ est la $\ZZ_\varphi$-extension de $K$ absolument galoisienne attachée à $\varphi$.
\end{Prop}

Ici $Z^\varphi$ est une $\ZZ_\ell^{d_\varphi}$-extension de $K$, où $d_\varphi=\deg\varphi$ est le degré du caractère $\varphi$, les caractères $\ell$-adiques irréductibles du groupe $\Delta$ n'étant absolument irréductibles que lorsque l'exposant de $\Delta$ divise $\ell-1$, auquel cas on a $\rho_\varphi^\tau(z)=\varphi(\tau)z$. Par exemple, un corps quadratique imaginaire possède ainsi exactement deux $\Zl$-extensions absolument galoisiennes: la cyclotomique et la prodiédrale.\medskip

\Preuve Ce résultat étant bien connu lorsque les caractères $\ell$-adiques irréductibles sont de degré 1 (cf. e.g. \cite{J29}, Prop. 9), expliquons rapidement comment l'obtenir sans cette hypothèse. Introduisons la pro-$\ell$-extension abélienne $\ell$-ramifiée maximale $M$ de $K$ et notons $H$ sa sous-extension maximale non-ramifiée. La théorie $\ell$-adique du corps de classes (cf. \cite{Gra2005} ou \cite{J31}) nous donne l'isomorphisme:
\smallskip

\centerline{$\Gal(M/H) \simeq \prod_\p\U_\p\,\R / \prod_{\p\nmid\ell}\U_\p\,\R \simeq \,\U_\ell/\E_\ell$,}\smallskip

\noindent où $\,\U_\ell=\prod_{\l|\ell}\U_\l$ est le groupe des unités semi-locales attaché aux places sauvages et $\,\E_\ell$ l'image dans $\,\U_\ell$ du $\ell$-adifié $\,\E=\Zl\otimes_\ZZ E$ du groupe des unités de $K$. Or, comme module noethérien sur $\Zl[\Delta]$, le numérateur $\,\U_\ell$ est le produit direct de son sous-module de torsion $\upmu^\ph_\ell=\prod_{\l|\ell}\upmu^\ph_\l$ et d'un module libre de rang 1. Et le dénominateur $\,\E_\ell$ est réel aux racines de l'unité près.
Prenant les $\varphi$-composantes, pour $\varphi$ irréductible imaginaire, nous avons donc le pseudo-isomorphisme:\smallskip

\centerline{$\Gal(M/K)^{e_\varphi} \sim \Gal(M/H)^{e_\varphi} \simeq (\U_\ell/\E_\ell)^{e_\varphi} \sim \ZZ_\varphi$.}\smallskip

\noindent Et $M$ contient bien une unique $\ZZ_\varphi$-extension, laquelle est, par construction, galoisienne sur $\QQ$. La proposition en résulte puisque toute $\ZZ_\ell^d$-extension de $K$ est $\ell$-ramifiée, donc contenue dans $M$.\medskip

Transportons alors aux pro-$\ell$-extensions abéliennes $N$ de $K$ qui sont absolument galoisiennes la décomposition canonique de $\Gal(N/K)$ comme produit de ses composantes isotypiques et écrivons ainsi $N=\prod_\varphi N^\varphi$ comme compositum des sous-extensions ainsi obtenues. Cela étant, nous avons:

\begin{Prop}\label{BP}
Désignons par $H$ la $\ell$-extension abélienne non-ramifiée maximale de $K$ (autrement dit son $\ell$-corps de classes de Hilbert), par $\wH$ son analogue logarithmique (i.e. la pro-$\ell$-extension abélienne logarithmiquement non-ramifiée maximale de $K$) et par $M_\bp$ la pro-$\ell$-extension de Bertrandias-Payan attachée à $K$.
Soit enfin $\varphi$ un caractère irréductible imaginaire de $\Delta$.\par
Avec les conventions ci-dessus, la sous-extension $\varphi$-isotypique $M_\bp^\varphi$ de $M_\bp$ de $K$ est le compositum de la $\ZZ_\varphi$-extension $Z^\varphi$ de sa sous-extension non-ramifiée $H^\varphi$ et de la $\ZZ_\varphi$-extension $Z^\varphi$ par et de la sous-extension $\varphi$-isotypique $H^\varphi$ de $H$ tout comme de celle $\wH^\varphi$ de $\wH$:\smallskip

\centerline{$M_\bp^\varphi = Z^\varphi H^\varphi = Z^\varphi \wH^\varphi$.}
\end{Prop}

\Preuve Rappelons que la pro-$\ell$-extension de Bertrandias-Payan $M_\bp$ de $K$ est le compositum des $\ell$-extensions cycliques de $K$ qui sont {\em localement} plongeables dans le compositum $Z$ des $\Zl$-extensions. En d'autres termes, $M_\bp$ est la pro-$\ell$-extension abélienne de $K$ fixée par le sous-groupe d'idèles $\prod_\p\upmu_\p\,\R$ (cf. \cite{J31}, Exemple 2.9). Et, comme on a $\,\U_\p=\upmu_\p$ pour $\p\nmid\ell$, elle coïncide avec la pro-$\ell$-extension abélienne $\ell$-ramifiée maximale $M$ dès lors que $\upmu_\ell/\upmu=(\prod_{\l|\ell}\upmu_\l)/\upmu$ est trivial.\par

La théorie $\ell$-adique du corps de classes (cf. \cite{J31}) nous donne alors les isomorphismes:\smallskip

\centerline{$\Gal(M_\bp/H) \simeq \prod_\p\U_\p\,\R / \prod_\p\upmu_\p\,\R \simeq \,\U_\ell/\upmu_\ell\E_\ell$ \quad \& \quad $\Gal(M_\bp/H) \simeq \prod_\p\wU_\p\,\R / \prod_\p\upmu_\p\,\R \simeq \,\wU_\ell/\upmu_\ell\wE_\ell$}\smallskip

\noindent puis, par restriction aux $\varphi$-composantes imaginaires:\smallskip

\centerline{$\Gal(M_\bp^\varphi/H^\varphi) \simeq (\U_\ell/\upmu_\ell)^{e_\varphi} \simeq \ZZ_\varphi$ \quad \& \quad $\Gal(M_\bp^\varphi/\wH^\varphi) \simeq (\wU_\ell/\upmu_\ell)^{e_\varphi} \simeq \ZZ_\varphi$.}\smallskip

\noindent Ainsi $M_\bp^\varphi/H^\varphi$ est une $\ZZ_\varphi$-extension qui contient $Z^\varphi H^\varphi/H^\varphi$. De même  $M_\bp^\varphi/\wH^\varphi$ et $Z^\varphi \wH^\varphi /\wH^\varphi$.
D'où la Proposition, deux $\ZZ_\varphi$-extensions emboîtées et de même base coïncidant nécessairement.

\newpage
\section{Extensions cyclotomiquement ramifiées}

La notion d'extension cyclotomiquement ramifiée a été introduite dans \cite{J29}. Rappelons en brièvement la définition ainsi que la description via le corps de classes (cf. \cite{J29}, Déf.1 \& Th.2):

\begin{DTh}
Une pro-$\ell$-extension abélienne $L$ d'un corps de nombres $K$ est dite cyclotomiquement ramifiée lorsqu'elle satisfait les deux propriétés suivantes:\par
$(i)$ Les places au-dessus de $\ell$ sont infiniment ramifiées dans $L/K$.\par
$(ii)$ L'extension induite $LK_\%/K_\%$ au-dessus de la $\Zl$-extension cyclotomique est non-ramifiée.\smallskip

\noindent Le groupe de normes associé à la pro-$\ell$-extension abélienne cyclotomiquement ramifiée maximale $\widehat{H}$ de $K$ est $\,\widehat{\U}\R$, où $\,\widehat{\U}=\U \cap \,\wU$ est formé des idèles unités au double sens ordinaire et logarithmique.
\end{DTh}

Notons, en effet, $K^\ab$ la pro-$\ell$-extension abélienne maximale de $K$. Par la théorie $\ell$-adique du corps de classes (cf. \cite{J31}), le groupe de décomposition de la place $\p$ dans $K^\ab/K_\%$ est donné par:\smallskip

\centerline{$D_\p(K^\ab/K_\%) \simeq (\R_\p\R \cap \wJ)/\R = \wU_\p\R/\R$;}

\noindent et le groupe d'inertie par:\par

\centerline{$I_\p(K^\ab/K_\%) \simeq (\U_\p\R \cap \wJ)/\R = \widehat{\U}_\p\R/\R$.}\medskip

Ce point précisé, il vient ici:

\begin{Prop}
Soit $K$ un corps abélien imaginaire et $\ell$ un nombre premier ne divisant pas son degré.
Pour chaque caractère $\ell$-adique irréductible imaginaire $\varphi\leqslant\chi_\ell^\ph$ du groupe $\Delta=\Gal(K/\QQ)$ la sous-extension $\varphi$-isotypique de $\widehat{H}$ coïncide celle de l'extension de Bertrandias-Payan $M_\bp$; et la $\ZZ_\varphi$-extension absolument galoisienne $Z^\varphi$  attachée à $\varphi$ est cyclotomiquement ramifiée.
\end{Prop}

\Preuve Observons d'abord que $\ZZ^\varphi$, en tant que $\ZZ_\ell^d$-extension, est infiniment ramifiée en l'une au moins $\l$ des places au-dessus de $\ell$. Étant galoisienne sur $\QQ$, elle est donc bien infiniment ramifiée en ses conjuguées, i.e. en {\em toutes} les places au-dessus de $\ell$. Sous l'hypothèse $\varphi \leq \chi_\ell^\ph$, i.e. pour $\varphi_{|_{\Delta_\ell}}=1$, $Z^\varphi$ provient d'une $\ZZ_\varphi$-extension du sous-corps de décomposition $K^{\Delta_\ell}$ de $\ell$ dans $K$. Ce n'est donc pas restreindre la généralité que supposer pour la démonstration $\ell$ complètement décomposé dans $K/\QQ$. Cela étant, nous avons alors $K_\l\simeq\Ql$ en chacune $\l$ des places au-dessus de $\ell$, donc:\smallskip

\centerline{$\R_\l=\U_\l\,\wU_\l = \upmu_\l (1+\ell_\l)_\ph^{\Zl} \ell_\l{}^{\Zl}$ \quad et \quad $\,\widehat{\U}_\l=\,\U_\l \cap\,\wU_\l=\upmu_\l$.}\smallskip

En particulier, il suit $M_\bp=\widehat{H}$, donc  $M_\bp^\varphi=\widehat{H}^\varphi$; puis le résultat par la Proposition \ref{BP}.
Et les extensions considérées prennent ainsi place dans le diagramme:\bigskip

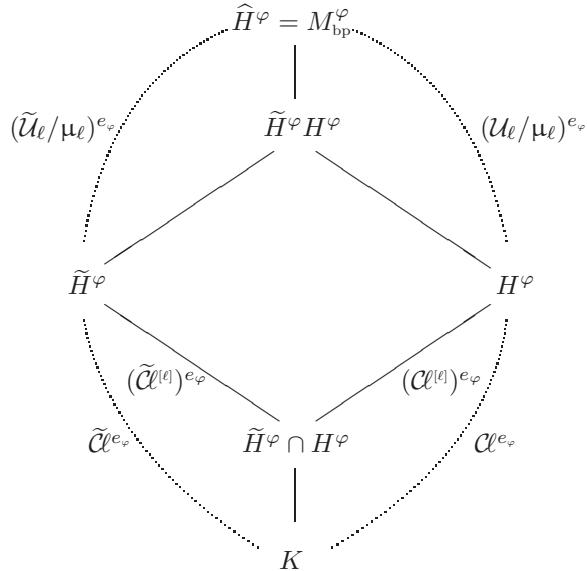
\begin{center}\label{Schema}
\setlength{\unitlength}{2pt}
\begin{picture}(106,115)(0,0)

\put(41,110){$\widehat{H}^\varphi=M_\bp^\varphi$}
\put(47,90){$\wH^\varphi H^\varphi$}
\put(10,60){$\wH^\varphi$}
\put(91,60){$H^\varphi$}
\put(43,30){$\wH^\varphi\cap H^\varphi$}
\put(50,8){$K$}

\bezier{60}(40,110)(18,100)(13,70)
\bezier{60}(64,110)(87,100)(93,70)
\bezier{60}(13,55)(18,30)(46,12)
\bezier{60}(93,55)(90,30)(60,12)

\put(53,107){\line(0,-1){10}}
\put(53,27){\line(0,-1){10}}
\put(49,87){\line(-3,-2){32}}
\put(57,87){\line(3,-2){33}}
\put(50,36){\line(-3,2){33}}
\put(57,36){\line(3,2){33}}

\put(-1,90){$(\wU_\ell/\upmu_\ell)^{e_{\si{\varphi}}}$}
\put(88,90){$(\U_\ell/\upmu_\ell)^{e_{\si{\varphi}}}$}
\put(14,30){$\wCl^{e_{\si{\varphi}}}$}
\put(87,30){$\Cl^{e_{\si{\varphi}}}$}

\put(21,42){$(\wCl^{\si{[\ell]}})^{e_\varphi}$}
\put(73,42){$(\Cl^{\si{[\ell]}})^{e_\varphi}$}

\end{picture}

\end{center}

\newpage
\section{Critère logarithmique dans le cas non semi-simple}

Revenons pour finir sur le le critère logarithmique sans l'hypothèse de semi-simplicité. Prenons toujours $\ell$ premier impair, mais considérons un corps abélien imaginaire arbitraire $L$.\smallskip

Supposons donc maintenant $\ell \mid [L:\QQ]$; notons $G$ le $\ell$-sous-groupe de Sylow de $\Gal(L/\QQ)$ et $\Delta$ son cofacteur; puis $K=L^G$ le sous-corps de $L$ fixé par $G$, qui est donc un corps abélien imaginaire de degré $d$ étranger à $\ell$ et de groupe de Galois $\Gal(K/\QQ)\simeq\Delta$.\smallskip

Tout comme pour $K$, les invariants arithmétiques attachés à $L$ étudiés ici sont canoniquement des $\Zl[\Delta]$-modules et nous pouvons leur appliquer les mêmes arguments liés à la décomposition semi-locale de l'algèbre $\Zl[\Delta]$ associée aux caractères $\ell$-adiques irréductibles du groupe $\Delta$. En particulier le Théorème \ref{CGold} s'applique {\em mutatis mutandis} à $L$. En revanche, le critère de trivialité de la $\varphi$-composante du $\ell$-groupe des classes logarithmiques donné par le Théorème \ref{CLog} ne peut se transposer directement, son sous-module sauvage n'étant plus $\ZZ_\varphi$-monogène. La formule des classes logarithmiques ambiges permet cependant de balayer cette difficulté. Il vient ainsi:

\begin{Th}\label{NSS}
Soient $L$ un corps abélien imaginaire, $\ell$ un nombre premier impair divisant son degré et $L_\%=\bigcup_n L_n$ la $\Zl$-extension cyclotomique de $L$. Notons $G$ le $\ell$-sous-groupe de Sylow de $\Gal(L/\QQ)$, puis $K$ le sous-corps de $L$ fixé par $G$ et $\Delta\simeq\Gal(K/\QQ)$ le cofacteur de $G$.\smallskip

Pour chaque caractère $\ell$-adique irréductible imaginaire $\varphi$, contenu dans l'induit $\chi^\ph_\ell=\Ind_{\Delta_\ell}^\Delta 1_{\Delta_\ell}$ du caractère unité du sous-groupe de décomposition de $\ell$ dans l'extension $K/\QQ$, la composante sauvage $\,\C_{L_\%}^{[\ell]\,e_\varphi}$ du pro-$\ell$-groupe $\,\C_{L_\%}^{\,e_\varphi}=\varprojlim \Cl_{L_n}^{\,e_\varphi}$ est un $\ZZ_\varphi$-module libre qui a pour dimension l'indice $d_\ell=\big(\Gal(L_\%/K)/D_\ell(\si{L_\%/K})\big)$ du sous-groupe de décomposition $D_\ell(L_\%/K)$ de $\ell$ dans $L_\%/K$. \smallskip

Par suite l'invariant lambda d'Iwasawa du $\Lambda_\varphi$-module $\,\C_L^{\,e_\varphi}$ vérifie l'inégalité $\lambda_L^\varphi \geqslant d_\ell$; et les assertions suivantes sont équivalentes:
\begin{itemize}\smallskip

\item[(i)] L'invariant structurel lambda du $\Lambda_\varphi$-module $\,\C_L^{\,e_\varphi} = \varprojlim \,\Cl_{L_n}^{\,e_\varphi}$ vaut: $\lambda_L^\varphi =d_\ell$.\smallskip

\item[(ii)] L'invariant structurel lambda du $\Lambda_\varphi$-module $\,\wC_L^{\,e_\varphi} = \varprojlim \,\wCl_{L_n}^{\,e_\varphi}$ vaut: $\wt\lambda_L^\varphi =0$.\smallskip

\item[(iii)] Le  $\Lambda_\varphi$-module $\,\wC_L^{\,e_\varphi}$ est trivial: $\,\wC_L^{\,e_\varphi} = 1$.\smallskip

\item[(iv)] La $\varphi$-composante du $\ell$-groupe des classes logarithmiques de $L$ est triviale: $\,\wCl_L^{\,e_\varphi} = 1$.\smallskip

\item[(v)] On a  $\,\wCl_K^{\,e_\varphi} = 1$ et la $\ell$-extension $L/K$ n' est ramifiée au sens logarithmique en aucune des places $p$ pour lesquelles $\varphi$ divise l'induit $\chi^\ph_p=\Ind_{\Delta_p}^\Delta 1_{\Delta_p}$ du caractère unité du sous-groupe de décomposition de $p$ dans $K/\QQ$.
\end{itemize}
\end{Th}

\Preuve Pour le premier point, toujours en l'absence d'unités imaginaires (y compris de racines de l'unité dans la $\varphi$-composante pour $\varphi\mid\chi_\ell^\ph$) et du fait que les places sauvages sont presque totalement ramifiées dans $L_\%/L$, on a banalement:\smallskip

\centerline{$\C_L^{\si{[\ell]\,e_\varphi}}\simeq\varprojlim \D_{L_n}^{\si{[\ell]\,e_\varphi}}\simeq\ZZ_\varphi[\Gal(L_\%/K)/D_\ell(L_\%/K)])\simeq\ZZ_\varphi^{d_\ell}$.}
\medskip

L'équivalence des propriétés $(i)$ à $(iv)$ s'établit alors tout comme pour le Théorème \ref{CGold}. Le seul point à vérifier est donc l'équivalence des assertions $(iv)$ et $(v)$. Pour voir cela, procédons par étapes en introduisant une suite d'extensions relativement cycliques $L^i/L^{i-\si{1}}$ allant de $L^{\si{0}}=K$ à $L^k=L$. Et appliquons aux $\varphi$-composantes la suite exacte des classes logarithmiques ambiges dans chacune des extensions cycliques $L^i/L^{i-\si{1}}$ pour $i=1,\cdots,k$. En l'absence d'unités logarithmiques imaginaires (y compris de racines de l'unité dans les $\varphi$-composantes par $\varphi\mid\chi_\ell^\ph$), il vient:\smallskip

\centerline{$|\,\wCl^{G_i\,e_\varphi}_{L^i}| = |\,\wCl^{\,e_\varphi}_{L^{i-\si{1}}}|\,\big( \wDl^{G_i\,e_\varphi}_{L^i}:\wDl^{\,e_\varphi}_{L^{i-\si{1}}} \big)$, avec $G_i=\Gal(L^i/L^{i-\si{1}})$.}\smallskip

\noindent D'où, de proche en proche, le résultat annoncé, les diviseurs logarithmiques ambiges étant exactement ceux construits sur les diviseurs logarithmiques étendus et ceux logarithmiquement ramifiés.\smallskip

\Remarque Le caractère abélien de $G$ n'est pas essentiel ici, les arguments développés dans la preuve du Théorème fonctionnant aussi bien dès lors que le corps $L$ est le compositum d'une $\ell$-extension (galoisienne) arbitraire de $\QQ$ et d'un corps abélien $K$ de degré étranger à $\ell$.

\newpage
\section{Lien avec les formules de Kida--Kuz'min}

Pour passer de $K$ à $L$ dans la preuve du Théorème \ref{NSS}, nous avons choisi d'utiliser la formule des classes ambiges pour les classes logarithmiques rappelée plus haut. Mais nous aurions pu tout aussi bien prendre appui sur les formules de transition de genre à la Kida--Kuz'min \cite{Ki1,Ki2,Ku}.\smallskip

Les invariants lambda attachés aux divers modules d'Iwasawa standard (non-ramifiés, décomposés partout, $\ell$-ramifiés, etc.) constituent, en effet, un analogue pour les corps surcirculaires, i.e. pour les $\Zl$-extensions cyclotomiques des corps de nombres, de ce qu'est le genre pour les corps de fonctions; ce qu'illustent notamment les formules de transition à la Riemann--Hurwitz (cf. \cite{J25}).\smallskip

Reprenons en effet, en les transposant dans le contexte de cette note, les résultats de \cite{J19} sur le genre des corps surcirculaires (Th. 23 et  Th. 27 dans la numérotation de la version Arxiv):

\begin{Th}\label{TG}
Soient $\ell$ un nombre premier impair, $\QQ_\%$ la $\Zl$-extension cyclotomique de $\QQ$ et $L$ une extension abélienne imaginaire de $\QQ$, composée directe d'une sous-extension imaginaire $K$ de degré étranger à $\ell$ et d'une $\ell$-extension réelle $F$. Notons $K'=K[\zeta_\ell]$ et $L'=L[\zeta_\ell]$ puis $\Delta'=\Gal(L'/K')$ et $\omega$ le caractère de Teichmüller.

Alors, pour chaque caractère $\ell$-adique irréductible imaginaire $\varphi$ de $\Delta'$, les invariants lambda d'Iwasawa respectifs de $K'$ et de $L'$ attachés aux $\varphi$-composantes sont liés par les identités:\smallskip

(i)\; $\lambda_L^\varphi -\langle\omega,\varphi\rangle\,\omega= [L_\%:K_\%]\,(\lambda_K^\varphi  -\langle\omega,\varphi\rangle\,\omega) + \sum_{p\ne\ell} \big(e_p(\si{L_\%/K_\%})-1\big) d_p(\si{K_\%/K})\frac{\langle\chi_p^\ph,\varphi\rangle}{\langle\varphi,\varphi\rangle}\,\varphi$, 
\smallskip

\noindent où la sommation porte sur les places $p$ modérément ramifiées dans $L_\%/K_\%)$; $e_p(\si{L_\%/K_\%})$ est l'indice de ramification correspondant; et $d_p(\si{K_\%/K})$ est l'indice de décomposition de $p$ dans $K_\%/K$.\smallskip

(ii)\; $\wt\lambda_L^\varphi -\langle\omega,\varphi\rangle\,\omega= [L_\%:K_\%]\,(\wt\lambda_K^\varphi  -\langle\omega,\varphi\rangle\,\omega) + \sum_{p} \big(\wt e_p(\si{L_\%/K_\%})-1\big) d_p(\si{K_\%/K})\frac{\langle\chi_p^\ph,\varphi\rangle}{\langle\varphi,\varphi\rangle}\,\varphi$, 
\smallskip 

\noindent où la sommation porte sur les places $p$ logarithmiquement ramifiées dans $L_\%/K_\%)$; $\wt e_p(\si{L_\%/K_\%})$ est l'indice de ramification correspondant; et $d_p(\si{K_\%/K})$ l'indice de décomposition de $p$ dans $K_\%/K$.
\end{Th}

\Preuve L'assertion $(i)$ n'est autre {\em mutatis mutandis} qu'une transcription du résultat donné par le Th. 23 de \cite{J19}.
L'assertion $(ii)$ appelle un petit commentaire: aux $p\ne\ell$, les indices de ramification, pris au sens ordinaire $e_p(\si{L_\%/K_\%})$ ou logarithmique $\wt e_p(\si{L_\%/K_\%})$, coïncident. En outre la montée dans la $\Zl$-extension cyclotomique ayant épuisé toute possibilité d'inertie (au double sens ordinaire et logarithmique), ils coïncident aussi avec le degré local $d_p(\si{L_\%/K_\%})$. Et le même résultat vaut pour l'indice logarithmique $\wt e_\ell(\si{L_\%/K_\%})=d_\ell(\si{L_\%/K_\%})$ en $p=\ell$. On peut donc remplacer $\wt e_p(\si{L_\%/K_\%})$ par $d_p(\si{L_\%/K_\%})$ dans $(ii)$ et retomber ainsi sur le résultat du Th. 27 de \cite{J19}.\smallskip

Dans les deux cas, le caractère $\varphi$ étant supposé imaginaire, il suit $\langle\chi_p^\ph,\varphi\rangle=0$ si les places au-dessus de $p$ sont invariantes par la conjugaison complexe $\bar\tau$. En d'autres termes, seules interviennent dans les formules $(i)$ et $(ii)$ les places $p$ décomposées par $\bar\tau$.\smallskip

Enfin, si $\varphi$ est représenté dans $\chi_\ell^\ph$, on a $\varphi\ne\omega$ et le terme correctif $\langle\omega,\varphi\rangle\,\omega$ est nul. Il vient donc, en cohérence avec le Théorème \ref{NSS}:

\begin{Cor}
Sous les hypothèses du Théorème \ref{NSS}, pour tout caractère irréductible imaginaire $\varphi$ de $\Delta=\Gal(L/H)$ représenté dans $\chi_\ell^\ph=\Ind_{\Delta_\ell}^\Delta 1_{\Delta_\ell}$, on a:\smallskip
\begin{itemize}
\item[(i)]\; $\lambda_L^\varphi = [L_\%:K_\%]\,\lambda_K^\varphi + \sum_{p\ne\ell} \big(e_p(\si{L_\%/K_\%})-1\big) d_p(\si{K_\%/K})\frac{\langle\chi_p^\ph,\varphi\rangle}{\langle\varphi,\varphi\rangle}\,\varphi$;
\smallskip
\item[(ii)]\, $\wt\lambda_L^\varphi = [L_\%:K_\%]\,\wt\lambda_K^\varphi   + \sum_{p} \big(\wt e_p(\si{L_\%/K_\%})-1\big) d_p(\si{K_\%/K})\frac{\langle\chi_p^\ph,\varphi\rangle}{\langle\varphi,\varphi\rangle}\,\varphi$.
\end{itemize}
\end{Cor}

En particulier, on retrouve ainsi l'équivalence donnée par le Théorème \ref{NSS}:
$$
\wt\lambda_L^\varphi = 0 \qquad \Leftrightarrow \qquad 
\left\{
\aligned
&\wt\lambda_K^\varphi = 0 \quad \textrm{et, pour tout } p \textrm { vérifiant } \varphi \mid \chi_p^\ph,\\
&L/K \text{ logarithmiquement non-ramifiée en }  p .
\endaligned
\right.
$$

Rappelons à cette occasion que {\em $L/K$ logarithmiquement non-ramifiée en $p$} signifie  {\em $L/K$ localement cyclotomique aux places au-dessus de $p$}, i.e. {\em $L/K$ localement contenue dans $K_\%/K$ en chacune des places $\p$ au-dessus de $p$}.

\newpage
\Remarques Pour les caractères irréductibles réels, la conjecture de Greenberg postule la trivialité des invariants $\lambda^\varphi$ (et, de façon équivalente, des invariants logarithmiques $\wi\lambda^\varphi$) pour chaque $\varphi$.\par

Pour les caractères irréductibles imaginaires, les formules obtenues montrent que la valeur minimale de l'invariant $\lambda^\varphi$ est déterminée explicitement par le schéma galoisien de la ramification; et ce minimum est réalisé en un $\varphi$ donné si et seulement si l'invariant logarithmique $\wi\lambda^\varphi$ est nul.\par

Enfin, les résultats précédents se transposent aisément aux composantes réelles du sous-module de $\Lambda$-torsion du module $\ell$-ramifié standard, via les identités de dualité discutées dans \cite{J43,J67,J40,J52} basées en toute généralité sur les théorèmes de réflexion de Gras \cite{Gra1998}.
\def\refname{\normalsize{\sc  Références}}
\addcontentsline{toc}{section}{Bibliographie}

{\footnotesize

\bigskip

\noindent{\sc Adresse:}
Univ. Bordeaux \& CNRS,\\
Institut de Mathématiques de Bordeaux,\\
351 Cours de la Libération,\\
F-33405 Talence cedex

\noindent{\sc Courriel:}
 \tt jean-francois.jaulent@math.u-bordeaux.fr
}

\end{document}